\documentclass[12pt]{article}
\usepackage{ifpdf}
\usepackage {graphics}
\usepackage{amsmath,mathptmx,amssymb,bm}
\newtheorem{theorem}{Theorem}

\newtheorem{definition}{Definition}

\begin{document}


\title{Symmetries and conservation laws of a fifth-order KdV equation with time-dependent coefficients and linear damping}


\author{R. de la Rosa${}^{a}$, M.L. Gandarias${}^{b}$, M.S. Bruz\'on${}^{c}$\\
 ${}^{a}$ Universidad de C\'adiz, Spain  (e-mail: rafael.delarosa@uca.es). \\
 ${}^{b}$ Universidad de C\'adiz, Spain (e-mail:  marialuz.gandarias@uca.es). \\
 ${}^{c}$ Universidad de C\'adiz, Spain  (e-mail:  m.bruzon@uca.es). \\  
}

\date{}
 
\maketitle

\begin{abstract}

A fifth-order KdV equation with time dependent coefficients and linear damping has been studied. Symmetry groups have several different applications in the context of nonlinear differential equations. For instance, they can be used to determine conservation laws. We obtain the symmetries of the model applying Lie's classical method. The choice of some arbitrary functions of the equation by the equivalence transformation enhances the study of Lie symmetries of the equation. We have determined the subclasses of the equation which are nonlinearly self-adjoint. This allow us to obtain conservation laws by using a theorem proved by Ibragimov which is based on the concept of adjoint equation for nonlinear differential equations.\\

\noindent \textit{Keywords}: Classical Symmetries; Equivalence Transformations; Partial Differential Equations; Conservation laws.\\
\noindent \textit{PACS}: 02.20.Hj; 02.20.Sv; 02.30.Em; 02.30.Jr; 47.10.ab.\\
\noindent \textit{MSC}: 35C07; 35Q53.
\end{abstract}



\maketitle

\section{Introduction}
\label{intro}

Lie symmetry groups play a fundamental role in the study of differential equations, particularly for nonlinear equations. Lie's theory provides a useful, powerful tool when analysing partial differential equations. Moreover, it has numerous well known applications, prominent among these are the obtaining of exact solutions of partial differential equations, directly or via similarity solutions \cite{gupta,wang}, clasify invariant equations, reduce the number of independent variables or determine conservation laws \cite{biswas,BGI,freire,gandarias,razborova}.

Over the last years, nonlinear equations with variable coefficients have become increasingly important due to these describe many nonlinear phenomena more realistically than equations with constant coefficients. For instance, the variable coefficient Gardner equation \cite{RGB15,krishnan}, the nonlinear Schr\"{o}dinger equation \cite{liu,senthil,zedan} or the variable coefficient reaction-diffusion equation \cite{BGR14,vaneeva}. The problem lies in the fact that the analysis of Lie symmetries of equations involving arbitrary functions seems rather difficult. Equivalence transformations allow to perform further study in a simpler way.

An equivalence transformation is a non-degenerate transformation acting on dependent and independent variables which maps the considered equation to another equation with the same differential structure, except maybe the form of the arbitrary functions. Most important advantage of equivalence transformations is that they permit an analysis for complete equivalent classes instead of considering individual equations. In particular, they allow us to reduce of a class to its subclass with fewer number of arbitrary functions.

The study of nonlinear evolution equations (NLEEs) is one of the most important areas of research in the field of nonlinear dynamics. There is a large number of nonlinear evolution equations which are being studied nowadays. Many of these equations are generalizations or combinations of the KdV equation, Boussinesq equation and many others.

The KdV equation is established as a simple mathematical model by which it can be described the physics of a shallow layer of fluid subject, for instance, it is an analytical model of tsunami generation by sub-marine landslides. On the other hand, the dynamics of shallow water waves can also be described by using a KdV equation. Generalizations of KdV equation which have more than one nonlinear term allow a better analysis of these phenomena.

The concept of conservation law, more specifically, of a conserved quantity, has its origin in physics. Conservation laws are very useful when it comes to analyse or know of physical properties of concrete systems, overall physical aspects of NLEEs. NLEEs would be difficult to understand without having deepened into the study of its conservation laws.

In this paper, we consider a family of fifth-order KdV equations with time-dependent coefficients and linear damping term of the form
\begin{equation}\label{eq}
\begin{array}{c} 
 u_t+A(t) u_{xxxxx}+B(t) u_{xxx}+ C(t) u u_{xxx}  \\ \qquad \quad +E(t) u u_x+ F(t) u_x u_{xx} + Q(t) u=0 ,
\end{array}\end{equation}
\noindent where  $A(t) \neq 0$, $B(t)$, $C(t) \neq 0$, $E(t)$, $F(t)$ and $Q(t)$ are arbitrary smooth functions of $t$. We perform an analysis of Lie symmetries of the equation. We obtain the continuous equivalence transformations of class (\ref{eq}). This allow us to reduce the number of arbitrary elements of class (\ref{eq}) which will result in a complete study of Lie symmetries of class (\ref{eq}). Moreover, we determine the subclasses of the equation which are quasi self-adjoint and nonlinear self-adjoint. From these concepts we construct conservation laws by using a theorem proved by Ibragimov.

\section{Equivalence Transformations}
\label{sec:1}
In this section, we search for equivalence transformation of class (\ref{eq}). An equivalence transformation of class (\ref{eq}) is a nondegenerate point transformation, $\left(t,x,u \right)$ to $\left(\tilde{t},\tilde{x},\tilde{u} \right)$ which remains the differential structure but with different arbitrary functions, $\tilde{A}(\tilde{t})$, $\tilde{B}(\tilde{t})$, $\tilde{C}(\tilde{t})$, $\tilde{E}(\tilde{t})$, $\tilde{F}(\tilde{t})$ and $\tilde{Q}(\tilde{t})$ from the original ones. In order to obtain the equivalence transformation for our equation we apply Lie's infinitesimal criterion \cite{ovsian}. In this case, we require not only the invariance of class (\ref{eq}) but also the invariance of the auxiliary system
\begin{equation}\label{aux}
\begin{array}{lll}
A_x & = & A_u=B_x=B_u=C_x=C_u= \\

& = & E_x=E_u=F_x=F_u=Q_x=Q_u=0
\end{array}
\end{equation}

\noindent in the augmented space $\left( t,x,u,A,B,C,E,F,Q  \right)$.\\

We consider the one-parameter group of equivalence transformations in $\left( t,x,u,A,B,C,E,F,Q  \right)$ given by
$$\begin{array}{rcl}\nonumber \tilde{t} & = & t+\epsilon \, \tau(t,x,u)+O(\epsilon^2),
\\\nonumber \tilde{x} & = & x+\epsilon \, \xi(t,x,u)+O(\epsilon ^2),\\\nonumber \tilde{u} & = & u+\epsilon \,
\eta(t,x,u)+O(\epsilon ^2),

\\\nonumber \tilde{A} & = & A+\epsilon \,
\omega^1( t,x,u,A,B,C,E,F,Q )+O(\epsilon ^2),
\\\nonumber \tilde{B} & = & B+\epsilon \,
\omega^2( t,x,u,A,B,C,E,F,Q  )+O(\epsilon ^2),
\\\nonumber \tilde{C} & = & C+\epsilon \,
\omega^3( t,x,u,A,B,C,E,F,Q  )+O(\epsilon ^2),
\\\nonumber \tilde{E} & = & E+\epsilon \,
\omega^4( t,x,u,A,B,C,E,F,Q  )+O(\epsilon ^2),
\\\nonumber \tilde{F} & = & F+\epsilon \,
\omega^5( t,x,u,A,B,C,E,F,Q  )+O(\epsilon ^2),
\\\nonumber \tilde{Q} & = & Q+\epsilon \,
\omega^6( t,x,u,A,B,C,E,F,Q  )+O(\epsilon ^2),
\end{array}$$ where  $\epsilon$ is the
group parameter. In this case, the vector field takes the following form\\

\begin{equation}\label{generator}
\begin{array}{rcl} Y & = &  \tau \partial_t + \xi \partial_x+\eta \partial_u+ \omega^1 \partial_A  + \omega^2 \partial_B  \\ \\ & &+ \omega^3 \partial_C + \omega^4 \partial_E + \omega^5 \partial_F + \omega^6 \partial_Q.  \end{array} \end{equation}

\noindent Invariance of system (\ref{eq})-(\ref{aux}) under a Lie group of point transformations with infinitesimal generator (\ref{generator}) leads to a system of determining equations. 

After having solved the determining system, omitting tedious calculations, we obtain that the equivalence generator  of class (\ref{eq}) is spanned by

$$\begin{array}{rcl}\nonumber Y_1 & = & x \partial_x+ \frac{1}{2}u \partial_u + 5 A \partial_A + 3 B \partial_B + \frac{5}{2} C \partial_C\\ \\

& & +  \frac{1}{2} E \partial_E+ \frac{5}{2} F \partial_F,\\ \\

\nonumber Y_2 & = & \partial_x,\\ \\

\nonumber Y_\alpha & = & \alpha \partial_t- \alpha_t A \partial_A- \alpha_t B \partial_B- \alpha_t C \partial_C \\ \\
& &- \alpha_t E \partial_E
- \alpha_t F \partial_F- \alpha_t Q \partial_Q,\\ \\

\nonumber Y_r & = &r u \partial_u - r C \partial_C- r E \partial_E- r F \partial_F- r_t \partial_Q,

\end{array}$$

\noindent where $\alpha=\alpha(t)$, $r=r(t)$ are arbitrary smooth functions with $\alpha_t \neq 0$. We obtain that the equivalence group of class (\ref{eq}) consists of the transformations\\
$$ \displaystyle \tilde{t}  =  \alpha(t), \quad \tilde{x}  =  (x+k_2)e^{k_1}, \quad \tilde{u}  =  e^{r(t)k_r +\frac{k_1}{2}}u, $$

\noindent where $k_1$, $k_2$ and $k_r$ are arbitrary constants, and the arbitrary functions are given by\\

$  \begin{array}{lllll}  \tilde{A}  =  \displaystyle \frac{e^{5 k_1}}{\alpha_t}A, & \quad &  \tilde{B}  =  \displaystyle \frac{e^{3 k_1}}{\alpha_t}B, \\ \\
 \tilde{C}   =   \displaystyle \frac{e^{\frac{5 k_1}{2} - r k_r}}{\alpha_t}C, & \quad & \tilde{E}  = \displaystyle \frac{e^{\frac{k_1}{2} - r k_r}}{\alpha_t}E, \\ \\
 \displaystyle  \tilde{F}  = \displaystyle \frac{e^{\frac{5 k_1}{2} - r k_r}}{\alpha_t}F,  & \quad & \tilde{Q}  = \displaystyle \frac{Q- r_t k_r}{\alpha_t}.
\end{array}$\\

\noindent As there are two arbitrary functions, $\alpha(t)$ and $r(t)$, in the group of transformations we can establish two of the arbitrary functions of class (\ref{eq}). Thus, we can set $\tilde{A}=\tilde{C}=1$ by the equivalence transformation

\begin{equation}\label{trans}
 \displaystyle \tilde{t}  =  \int e^{5 k_1} A dt, \,  \tilde{x}  =  (x+k_2)e^{k_1}, \, \tilde{u}  =  \frac{e^{-2 k_1}C}{A}u, 
\end{equation}

\noindent which changes equation (\ref{eq}) into the equation \\ 

$\begin{array}{c} 
 \tilde{u}_{\tilde{t}}+ \tilde{u}_{\tilde{x}\tilde{x}\tilde{x}\tilde{x}\tilde{x}}+\tilde{B}(\tilde{t}) \tilde{u}_{\tilde{x}\tilde{x}\tilde{x}}+ \tilde{u} 
 \tilde{u}_{\tilde{x}\tilde{x}\tilde{x}}  \\ \qquad \quad +\tilde{E}(\tilde{t}) \tilde{u} \tilde{u}_{\tilde{x}}+ \tilde{F}(\tilde{t}) \tilde{u}_{\tilde{x}} \tilde{u}_{\tilde{x}\tilde{x}} + \tilde{Q}(\tilde{t}) \tilde{u}=0 ,
\end{array}$\\

\noindent where the new arbitrary functions are related with the old ones by\\

$  \begin{array}{lllll}   \tilde{B}  =  \displaystyle \frac{e^{-2 k_1}B}{A}, & \quad & \tilde{E}  = \displaystyle \frac{e^{-2 k_1}E}{C}, \\ \\
\tilde{F}  = \displaystyle \frac{F}{C},  & \quad & \tilde{Q}  = \displaystyle \frac{e^{-5 k_1}\left(Q+ \left( Log\left( \frac{A}{C}\right) \right)_t  \right)}{A}.
\end{array}$\\

\noindent Thus, we can restrict without loss of generality our study to the class 

\begin{equation}\label{eqTR}
\begin{array}{c} 
 u_t+u_{xxxxx}+B(t) u_{xxx}+ u u_{xxx}  \\ \qquad \quad +E(t) u u_x+ F(t) u_x u_{xx} + Q(t) u=0 ,
\end{array}\end{equation}

\noindent due to symmetries and conservation laws obtained for class (\ref{eqTR}) can be extended to class (\ref{eq}) using transformation (\ref{trans}). 

\section{Lie classical symmetries of class (\ref{eqTR})}\label{symm}

Lie's classical method is based on the determination of the symmetry group of a differential equation, i.e., the largest group of transformations acting on dependents and independents variables of the equation so that transforms solutions of the equation into other solutions.\\

\noindent We consider the one-parameter Lie group of infinitesimal transformations in
$(t,x,u)$ given by
$$\begin{array}{lcr}\nonumber t^* & = & t+\epsilon \tau(t,x,u)+O(\epsilon ^2),
\\\nonumber x^* & = & x+\epsilon \xi(t,x,u)+O(\epsilon^2),\\\nonumber u^* & = & u+\epsilon
\eta(t,x,u)+O(\epsilon ^2),\end{array}$$ where  $\epsilon$ is the
group parameter.\\

\noindent The infinitesimal operator of the local Lie group of point transformations which are admitted by equation (\ref{eqTR}) is given by

\begin{equation}\label{vect}{\bf v}=
\tau(t,x,u)\partial_t+\xi(t,x,u) \partial_x+\eta(t,x,u)\partial_u.\end{equation} 

\noindent The fifth prolongation of the vector field (\ref{vect}) has the following form 

\begin{equation}\label{prolong}  \begin{array}{c} pr^{(5)}{\bf v}  ={\bf v} + \zeta^t \partial u_t+ \zeta^x \partial u_x +\zeta^{xx} \partial u_{xx} \\ \qquad +\zeta^{xxx} \partial u_{xxx}+\zeta^{xxxxx} \partial u_{xxxxx}, \end{array} \end{equation} 

\noindent where $$\zeta^J(t,x,u^{(5)})=D_J(\eta-\tau u_t-\xi u_x)+\tau
u_{Jt} +\xi u_{Jx}$$
with $J=(j_1,\ldots,j_k)$,  $1\leq j_k\leq 2$ and $1\leq
k\leq 5$, and $u^{(5)}$ denote the sets of partial derivatives up to fifth order \cite{olver}.
Invariance criterion implies that the 5th-prolongation of the vector field {\bf v} (\ref{prolong}) acting on the equation (\ref{eqTR}) must be zero when the equation holds. This leads to a set of 17 determining equations. By simplifying the system we obtain that the infinitesimals are given by
$$\tau= k_2 t+ k_3, \quad \xi= \frac{k_2}{5}x+ \delta, \quad \eta= \rho - \frac{2 k_2}{5}u, $$
where $\delta=\delta(t)$, $\rho=\rho(t,x)$, $B=B(t)$, $E=E(t)$, $F=F(t)$ and $Q=Q(t)$ must satisfy the following conditions:\\

\begin{equation}\label{sis}\begin{array}{r}

5(k_2 t+k_3)B_t+2 k_2 B +5 \rho=0,\\

\rho_x F=0,\\

\displaystyle (k_2 t+k_3)E_t+\frac{2 k_2}{5} E=0,\\

\rho_{xx} F+\rho E- \delta_t=0,\\

(k_2 t+k_3)F_t=0,\\

\rho Q+ \rho_{xxx}B+\rho_{xxxxx}+\rho_t=0, \\

(k_2 t+k_3)Q_t+k_2 Q+\rho_x E+\rho_{xxx}=0.\\

\end{array} \end{equation}

\noindent From determining system (\ref{sis}) if the functions $B$, $E$, $F$, $Q$, $\delta$ and $\rho$ are arbitrary we obtain ${\bf v}= \partial_x$. For the following cases we obtain new symmetries:

\subsection{Case 1: $F = 0$}
In this case, we get the following symmetries
$$\tau= k_2 t+ k_3, \quad \xi= \frac{k_2}{5}x+ \delta, \quad \eta= \rho - \frac{2k_2}{5}u, $$
where $\delta=\delta(t)$, $\rho=\rho(t,x)$, $B=B(t)$, $E=E(t)$ and $Q=Q(t)$ must satisfy the following conditions:\\
\begin{eqnarray}{}
\label{eq1c2} 5(k_2 t+k_3)B_t+2 k_2 B +5 \rho &=&0,\\
\label{eq2c2} \displaystyle (k_2 t+k_3)E_t+\frac{2 k_2}{5} E &=&0,\\
\label{eq3c2}\rho E- \delta_t &=&0,\\
\label{eq4c2}\rho Q+ \rho_{xxx}B+\rho_{xxxxx}+\rho_t &=&0,\\
\label{eq5c2} (k_2 t+k_3)Q_t+k_2 Q+\rho_x E+\rho_{xxx}&=&0.
\end{eqnarray}

\subsection{Case 2: $F = \mbox{constant} \neq 0 $}
Now, we obtain the following symmetries
$$\tau= k_2 t+ k_3, \quad \xi= \frac{k_2}{5}x+ \delta, \quad \eta= \sigma - \frac{2k_2}{5}u, $$
where $\delta=\delta(t)$, $\sigma=\sigma(t)$, $B=B(t)$, $E=E(t)$, $F=F(t)$ and $Q=Q(t)$ are related by the following conditions:\\
\begin{eqnarray}{}
\label{eq1c1} 5(k_2 t+k_3)B_t+2 k_2 B +5 \sigma &=&0,\\
\label{eq2c1} \displaystyle (k_2 t+k_3)E_t+\frac{2 k_2}{5} E &=&0,\\
\label{eq3c1}\sigma E- \delta_t &=&0,\\
\label{eq4c1}\sigma Q+\sigma_t &=&0,\\
\label{eq5c1} (k_2 t+k_3)Q_t+k_2 Q&=&0.
\end{eqnarray}

\noindent In the previous cases, $k_2$ and $k_3$ are arbitrary constants.

\section{Formal Lagragian, adjoint equation and nonlinearly self-adjointness}

In this section we proceed to give some prior concepts necessary to construct conservation laws. To achieve this objective, we use a theorem on conservation laws proved by Ibragimov \cite{Ibra:07}. The theorem is valid for any system of differential equations
where the number of equations is equal to the number of dependent variables, and does not require existence of a classical Lagrangian. This theorem is based on the concept of adjoint equation for nonlinear equations which is defined as

\begin{definition}\label{def2} Consider an sth-order partial differential equation
 \begin{equation}
 \label{fa}
 F({\rm x},u,u_{(1)}, \ldots,u_{(s)}) =0
 \end{equation}
 with independent variables ${\rm x}=(x^1,\ldots,x^n)$ and
 a dependent variable $u,$ where
$u_{(1)}=\{u_i\},$ $u_{(2)}=\{u_{ij}\},\ldots\,$
 denote the sets of the
partial derivatives  of the first, second, etc. orders,
$u_i=\partial u/\partial x^i$, $u_{ij}=\partial^2 u/\partial
x^i\partial x^j.$ The formal Lagrangian is defined as
\begin{equation}\label{lag}{\cal L}=v\,F\left({\rm x},u, u_{(1)}, \ldots,u_{(s)}\right), \end{equation} where
$v=v({\rm x})$ is a new dependent variable.  The adjoint equation to
{\rm(\ref{fa})} is
 \begin{equation}
 \label{faadj}
 F^{*}({\rm x},u,v,u_{(1)}, v_{(1)},\ldots,u_{(s)},v_{(s)})
 =0, \end{equation}
 with
 $$F^{*}({\rm x},u,v,u_{(1)}, v_{(1)},\ldots,u_{(s)},v_{(s)})
 =\frac{\delta(v\,F)}{\delta u},$$
where
$$ \frac{\delta }{\delta u}=\frac{\partial}{\partial u}
 +\displaystyle\sum_{s=1}^{\infty} (-1)^s D_{i_1}\cdots
D_{i_s}\frac{\partial}{\partial u_{i_1\cdots i_s}}$$ denotes the
variational derivatives (the Euler-Lagrange
 operator).
 Here
  $$
 D_i=\displaystyle\frac{\partial}{\partial
x^i}+u_i\frac{\partial}{\partial u}+ u_{ij}\frac{\partial}{\partial
u_j}+\cdots
 $$
 are the total differentiations.
\end{definition}

\begin{theorem}
The adjoint equation to class (\ref{eqTR}) is
\begin{equation}\label{adjointTR}
\begin{array}{c} 
v Q + (u_x v_{xx}+u_{xx}v_x) F-u v_x E-v_{xxx}B \\ \quad -v_{xxxxx} -u v_{xxx}-3 u_x v_{xx}-3u_{xx} v_x-v_t=0 
\end{array}
\end{equation}
\end{theorem}

\noindent Now, we use the following definition given in \cite{Ibran:11}

\begin{definition} Equation {\rm (\ref{fa})} is said to be {\bf nonlinearly self-adjoint} if the
 equation obtained from the adjoint equation {\rm (\ref{faadj})}
 by the substitution \begin{equation}\label{newv}v=\varphi({\rm x},u),\end{equation}
 with a certain function $\varphi({\rm x},u)\neq 0$
 is identical with the original equation {\rm(\ref{fa})}, i.e.,
 \begin{equation}\label{condadj}
 F^*\,{|}_{v=\varphi}=\lambda ({\rm x},u,...) F,\end{equation}
\noindent for some differential function $ \lambda= \lambda({\rm x},u,...)$. If $\varphi=u$  or $\varphi=\varphi(u)$ and $\varphi'(u)\neq 0$, equation {\rm (\ref{fa})} is said {\bf self-adjoint} or {\bf quasi self-adjoint}, respectively.\end{definition}

\noindent Given equation (\ref{eqTR}) we apply Definition 2. Taking into account expression (\ref{adjointTR}) and using (\ref{newv}) and its derivatives, equation (\ref{condadj}) can be written as

\begin{equation}\label{eqlambda}\begin{array}{l} -u Q \lambda+\varphi Q-\varphi_{x} u E-\varphi_{xxx} B-\varphi_{x x x}\,u-\varphi_{x x x x x}\\ -\varphi_{t} -\left(\lambda+\varphi_{u}\right) u_t- \varphi_{uuuuu} u_x^5 - 5\varphi_{uuuux} u_x^4  \\+ \left( \varphi_{uu}F-\varphi_{uuu}B-\varphi_{uuu}u-10\varphi_{uuuxx}-3\varphi_{uu}\right) u_x^3 \\ + \left( 2 \varphi_{u x} F-3 \varphi_{u u x} B-3 \varphi_{u u x} u-10 \varphi_{u u x x x} \right. \\ \left. -6 \varphi_{u x} \right) u_x^2 - \left( u E \lambda-\varphi_{x x}\,F+\varphi_{u} u E+3 \varphi_{u xx} B \right. \\ \left. +3 \varphi_{u x x} u +3 \varphi_{x x} +5 \varphi_{u x x x x} \right) u_x -15\varphi_{uux}u_{xx}^2 \\ - \left(B \lambda+u \lambda  +\varphi_{u} B +\varphi_{u} u  +10 \varphi_{u x x}\right) u_{xxx} \\ + \left( \varphi_{x}F-3 \varphi_{u x} B -3 \varphi_{u x} u-3 \varphi_{x}-10 \varphi_{u x x x} \right) u_{xx} \\ -10 \varphi_{uu}u_{xx}u_{xxx} -10 \varphi_{uuuu}u_{x}^3u_{xx}-10 \varphi_{uuu}u_{x}^2 u_{xxx} \\ -30 \varphi_{uuux}u_{x}^2 u_{xx}  -\left( \lambda+\varphi_{u} \right) u_{xxxxx} -5 \varphi_{ux} u_{xxxx}\\ -5 \varphi_{uu} u_x u_{xxxx} -15 \varphi_{uuu} u_x u_{xx}^2-20 \varphi_{uux} u_x u_{xxx} \\  - \left(F \lambda-2 \varphi_{u} F+3 \varphi_{u u} B+3 \varphi_{u u } u+30 \varphi_{u u x x} \right.\\ \left. +6 \varphi_{u} \right) u_x u_{xx}=0  .\end{array} \end{equation}

\noindent From the $u_t$, $u_x$, $u_{xx}$,... coefficients we obtain

\begin{theorem} Equation (\ref{eqTR}) with $B(t)$, $E(t)$, $F(t)$ and $Q(t)$ arbitrary functions is nonlinearly self-adjoint and
\begin{equation}\label{hs}\varphi=\alpha(t)u+\beta(t,x)\end{equation}
in the following cases:\\
\begin{itemize}
\item[$\bullet$] If $F(t)=2$, we obtain
$$\displaystyle \alpha(t)=c_1 e^{2\int Q(t) dt}, \quad  \beta(t)=c_2 e^{\int Q(t) dt}. $$
\item[$\bullet$] If $F(t)=3$, we obtain
$$\displaystyle \alpha(t)=0, \quad  \beta=\beta(t,x),$$
\noindent where $\beta(t,x)$ must verify the following equation
\begin{equation}\label{condicion} \beta Q+ \beta_x \, B \,E - \beta_{xxxxx}-\beta_t=0. \end{equation}
\item[$\bullet$] If $F(t)\neq 2, 3$, we obtain
$$\displaystyle \alpha(t)=0, \quad  \beta(t)=c_2 e^{\int Q(t) dt}. $$
\end{itemize}
Here, $c_1$ and $c_2$ are arbitrary constants.
\end{theorem}

\section{Conservation laws}

Conservation laws appear in many of physical, chemical and mechanical processes, such laws enable us to solve problems in which certain physical properties do not change in the course of time within an isolated physical system. The importance of conservation laws also embraces mathematics, for instance, the integrability of a partial differential equation (PDE) is strongly related with the existence of conservation laws. Furthermore, they can be used to obtain exact solutions of a PDE.\\

\noindent Given equation (\ref{eqTR}), a conservation law is a space-time divergence such that
\begin{equation}\label{cl}
\displaystyle{D_t C^1(t,x,u,u_x,u_t,...)+D_x C^2(t,x,u,u_x,u_t,...)=0.}
\end{equation}

\noindent In order to obtain conservation laws we use the following theorem given in \cite{Ibra:07}

\begin{theorem} Any Lie point, Lie-B\"{a}cklund or non-local symmetry
 \begin{equation}
 \label{gLC}
 {\bf v}=\xi^i({\rm x},u,u_{(1)},\ldots)\frac{\partial}{\partial x^i}
 +\eta({\rm x},u,u_{(1)},\ldots)\frac{\partial}{\partial u},
 \end{equation} of  equation {\rm (\ref{fa})}
 provides a conservation law $\,D_i(C^i)=0\, $ for the simultaneous system {\rm (\ref{fa})}, {\rm
 (\ref{faadj})}. The conserved vector is given by
  \begin{equation}
 \label{e19}
\begin{array}{ll}
 C^1&=\tau{\cal L}+ W\left[\displaystyle\frac{\partial{\cal L}}{\partial u_t}\right],\\[3.5ex]
 
 C^2 &= \xi{\cal L}+ W\left[\displaystyle \frac{\partial{\cal L}}{\partial u_x}-D_x\left( \frac{\partial {\cal L}}{\partial u_{xx}} \right) +D_x^2 \left(\frac{\partial {\cal L}}{\partial u_{xxx}}\right)\right. \\[3.5ex] &
 \left. \displaystyle  +D_x^4 \left(\frac{\partial {\cal L}}{\partial u_{xxxxx}}\right) \right] 
 + D_x \left( W \right) \left[\displaystyle \frac{\partial{\cal L}}{\partial u_{xx}}  \right.\\[3.5ex] &
\left. \displaystyle -D_x\left( \frac{\partial {\cal L}}{\partial u_{xxx}} \right) -D_x^3 \left( \frac{\partial {\cal L}}{\partial u_{xxxxx}} \right) \right]  \\[3.5ex] & 
+D_x^2\left(W\right) \left[ \displaystyle \frac{\partial {\cal L}}{\partial u_{xxx}} + D_x^2 \left(  \frac{\partial {\cal L}}{\partial u_{xxxxx}} \right) \right] \\[3.5ex] & 
\displaystyle -D_x^3 \left( W \right) \left[ D_x \left( \frac{\partial {\cal L}}{\partial u_{xxxxx}} \right) \right] +D_x^4 \left( W \right) \frac{\partial {\cal L}}{\partial u_{xxxxx}},
 \end{array}
 \end{equation}
 where ${\cal L}$ is given by (\ref{lag}) and $W$ is defined as follows:
 \begin{equation}
 \label{la}
 W = \eta - \xi^j u_j.
 \end{equation}
 \end{theorem}

\noindent In order to apply Theorem 3 to our equation we perform the following change of notation:
$\left(  \xi^1, \xi^2  \right)= \left(  \tau , \xi    \right).$ \\

\noindent Now, we proceed to construct conservation laws for the cases considered in Section \ref{symm}. We omit again the case in which the functions are arbitrary due to this leads to a trivial conservation law.

\subsection{Case 1.} In this case, we obtain conservation laws for $F(t)=0$ and $F(t)=\mbox{constant} \neq 2,3$. For these cases, the equation admits the following generator
$${\bf v}= \left(k_2 t+ k_3\right)\partial_t +\left( \frac{k_2}{5}x+ \delta \right) \partial_x+\left(\rho - \frac{2k_2}{5}u \right)\partial_u, $$
where $\delta=\delta(t)$ and $\rho=\rho(t,x)$ satisfy the conditions (\ref{eq1c2})-(\ref{eq5c2}) for $F(t)=0$, $\delta=\delta(t)$ and $\rho=\sigma(t)$ verify the conditions (\ref{eq1c1})-(\ref{eq5c1}) for $F(t) =\mbox{constant} \neq 2,3$. From Theorem 2, for both cases, we have
$$\varphi= c_2 e^{\int Q(t) dt}$$
Thus, we obtain the conservation law (\ref{cl}) with the conserved vector
$$ \begin{array}{lll} C^1 & = & \displaystyle c_2 e^H \left( \rho- \frac{k_2}{5}u+ (k_2 t+ k_3) Q u   \right), \\ \\

C^2 & = & \displaystyle \frac{c_2}{2} \left( k_2 t+k_3 \right)  e^H Q \left( E u^2 +F u_x^2-2 B u_{xx}-2 u_{xxxx} \right. \\ \\
& & \displaystyle \left. -u_x^2 -2 u u_{xx}  \right) -\frac{c_2 k_2}{10} e^H  \left(  2 B u_{xx}+F u_x^2+2 u_{xxxx} \right) \\ \\
& & +2 u u_{xx}-u_x^2 +c_2 e^H \left( \rho_{xx}B-\rho_x u_x+ \rho_{xx}u+\rho_{xxxx}  \right),\end{array}
$$
where $H(t)=\int Q(t) dt$. Remember that for case $F(t)=\mbox{constant} \neq 2,3$  the terms in which appear $\rho_x$, $\rho_{xx}$,... are eliminated.

\subsection{Case 2.} Now, we obtain conservation laws of equation (\ref{eqTR}) for $F(t)=2$ and $F(t)=3$. For these cases, the equation admits the following generator\\
$${\bf v}= \left(k_2 t+ k_3\right)\partial_t +\left( \frac{k_2}{5}x+ \delta \right) \partial_x+\left(\sigma - \frac{2k_2}{5}u \right)\partial_u, $$\\
where $\delta=\delta(t)$ and $\sigma=\sigma(t)$ must verify the conditions (\ref{eq1c1})-(\ref{eq5c1}).

\subsubsection{$F(t)=2$} We focus on case $F(t)=2$. From Theorem 2, we have
$$\varphi= c_1 e^{2 \int Q(t) dt}u+ c_2 e^{\int Q(t) dt}.$$
By using Theorem 3, we obtain the conserved vector
$$ \begin{array}{lll} C^1 & = & \displaystyle \left( c_1 e^H u +c_2 \right)(k_2 t+ k_3) e^H Q u \\ \\ & & \displaystyle + c_1 e^{2 H} u \left( \sigma - \frac{3 k_2}{10}u  \right)+ c_2 e^H \left( \sigma-\frac{k2}{5}u  \right), \\ \\

C^2 & = & \displaystyle \left( k_2 t+k_3 \right)  e^H Q \left( \left( 2 u u_{xxxx}-2 u_x u_{xxx}+u_{xx}^2  \right. \right. \\ \\
& & \displaystyle  \left. +2 u^2 u_{xx} -B u_x^2 +2 B u u_{xx} +\frac{2}{3}E u^3\right) c_1 e^H \\ \\
& & \displaystyle +\left. \left( \frac{1}{2}E u^2+B u_{xx}+u_{xxxx}+u u_{xx}+u_x^2 \right) c_2 \right) \\ \\

& & \displaystyle+\frac{c_1}{2}e^H \sigma \left( E u^2+ 2 B u_{xx}+2 u_{xxxx}+2 u u_{xx}+ u_x^2 \right) \\ \\

& &\displaystyle -\frac{3}{10} \left(  \frac{2}{3}E u^3+ 2 B u u_{xx} -B u_x^2+2 u u_{xxxx} +u_{xx}^2 \right. \\ \\

& & \displaystyle\left. -2 u_x u_{xxx}+2 u^2 u_{xx} \right) -\frac{3}{10} \left( E u^2 + 2 B u_{xx}+ u_x^2 \right. \\ \\

& & \left. + 2 u_{xxxx}+2 u u_{xx} \right) ,\end{array}
$$
where $H(t)=\int Q(t) dt$.

\subsubsection{$F(t)=3$} Now, we consider case $F(t)=3$. From Theorem 2, we have
$$\varphi= \beta (x,t),$$
where $\beta$ must verify (\ref{condicion}). By using Theorem 3, we obtain the conserved density
$$ \begin{array}{lll} C^1 & = & \displaystyle \left(k_2 t+k_3 \right) u \beta_t+ \left( \frac{k_2 x}{5}+\delta \right) \beta_x+ \left(\sigma- \frac{k_2 u}{5} \right) \beta, \\ \\
\end{array}$$

\noindent We do not show the conserved flux $C^2$ due to its length but we can state that the conserved vector obtained verifies (\ref{cl}).

\section{Conclusions}

In this paper, Lie group analysis has been successfully applied to study the symmetries of a fifth-order KdV equation with time dependent coefficients and linear damping. We have obtained the infinitesimal generators of the equivalence transformation of the equation. It is shown that the equivalence group obtained is an infinite dimensional group. Moreover, the usage of equivalence groups allow an exhaustive study of the equation transform class (\ref{eq}) into subclass (\ref{eqTR}) of simpler structure. We have obtained classical symmetries of class (\ref{eqTR}) involving the different arbitrary functions. Nonlinear self-adjointness of equation (\ref{eqTR}) has been carried out. This allow us to obtain new conservation laws by using Ibragimov's method. 

\section*{Acknowledgements}
We would like to thank the Editor and Referees for their timely and valuable comments and suggestions. The authors acknowledge the financial support from Junta de Andalucia group FQM-201. The first author express his sincerest gratitude to the Universidad Polit\'{e}cnica de Cartagena for support him. The second and third authors also acknowledge the support of DGICYT project MTM2009-11875 with the participation of FEDER.



\end{document}